\documentclass[a4paper]{amsart} 
\usepackage{amssymb} 
\newtheorem{theorem}{Theorem} 
\newtheorem{corollary}{Corollary}

\theoremstyle{definition} 
\newtheorem*{definition}{Definition} 
\newtheorem{remark}{Remark}

\begin{document} 
\title[Unconditionally $\tau$-Closed and $\tau$-Algebraic Sets in Groups]
{Unconditionally $\tau$-Closed\\ 
and $\tau$-Algebraic Sets in Groups}
\author{Ol'ga V. Sipacheva} 
\begin{abstract} 
Families of unconditionally $\tau$-closed and $\tau$-algebraic 
sets in a group are defined, 
which 
are natural generalizations of 
unconditionally closed and algebraic sets defined by Markov. 
A sufficient condition for the coincidence  of these families 
is found. 
In particular, it is proved that these families coincide
in any group of cardinality at most $\tau$. This result generalizes 
both Markov's theorem on the coincidence of unconditionally closed and 
algebraic sets in a countable group (as is known, they may be different 
in an uncountable group) and Podewski's theorem on the topologizablity 
of any \emph{ungebunden} group.
\end{abstract} 

\thanks{This work was financially 
supported by the Russian Foundation for Basic Research (project 
no.~06-01-00764).} 

\address {Department of General Topology and 
Geometry\\ Mechanics and Mathematics Faculty\\ Moscow State 
University\\ Leninskie Gory\\ Moscow, 119992 Russia}


\email
{o-sipa@yandex.ru}

\maketitle 
 
Markov~\cite{Markov1945} 
called a subset $A$ of a group $G$ 
\emph{unconditionally closed} in $G$ if it is closed in 
each 
Hausdorff group topology on $G$.  
Thus, a group is \emph{topologizable} 
(i.e.,  it admits a nondiscrete Hausdorff group topology) if and only if 
the complement to the identity (or any other) element in this group is 
not 
unconditionally closed. 

A natural example of unconditionally closed sets is 
the solution sets of equations in $G$, as well 
as their finite unions and arbitrary 
intersections. Markov called such sets algebraic. The precise 
definition is as follows. 

\begin{definition}[Markov~\cite{Markov1945}] 
A subset $A$ of a group $G$ with identity element $1$ is said to be 
\emph{elementary algebraic} in $G$ if there exists a word $w= w(x)$ in the 
alphabet $G\cup \{x^{\pm1}\}$ ($x$ is a variable) such that 
$A =\{x\in G: w(x) = 1\}$. 
Finite unions of elementary algebraic sets are \emph{additively 
algebraic} sets. Arbitrary 
intersections of additively algebraic sets are called \emph{algebraic}. 
\end{definition} 

Note that $G\setminus \{1\}$ is algebraic if and only 
if it is additively algebraic. (Indeed, if 
$G\setminus \{1\}= \bigcap A_\alpha$, where each $A_\alpha$ is a subset of 
$G$, then $A_\alpha = G\setminus \{1\}$ for some $\alpha$ 
and all the other $A_\beta$ coincide with $G$.)

An expression of the form $w(x)=1$, where $w(x)$ is an element of 
the free product $G*\langle x\rangle$,  i.e., a word 
in the alphabet $G\cup \{x^{\pm1}\}$ ($x$ is treated as a variable), 
is called an \emph{equation} 
in $G$. A \emph{solution} to this equation is any $a\in G$ such that 
substituting it for $x$ into $w(x)$ yields~1, i.e., such that $w(a)= 1$. 
An \emph{inequation} is an expression of the form $w(x)\ne1$,  
and its solution is any $a\in G$ for which $w(a)\ne 1$. 
The algebraic sets in $G$ are the solution sets of arbitrary 
conjunctions of finite disjunctions of equations in $G$. 

Thus, a natural necessary condition for a group $G$ to be 
topologizable is that the complement to the identity 
(or, equivalently, an arbitrary) element 
must not be (additively) algebraic in $G$ (this 
means that a finite system 
of inequations in $G$ cannot have precisely one solution). 

It had long been unknown whether this condition is also sufficient,  
i.e., whether the complement to the identity element in a group is 
unconditionally closed only if it is algebraic. 
In his 1945 paper~\cite{Markov1945}, Markov posed the problem of 
whether any unconditionally closed set is algebraic. 
In~\cite{Markov1946} (see also~\cite{Markov1944}), he solved this 
problem for countable groups by proving that any unconditionally 
closed set in a countable group is algebraic. 
This author recently proved that the answer is also positive 
for some 
subgroups of direct products of countable groups (in 
particular, for all Abelian groups)~\cite{sipa}. However, in the 
general case, the answer is negative (\cite{Hesse}; see also \cite{Sipa}).
 
Thus, the necessary topologizability condition stated above is not 
sufficient. In~1977, Podewski suggested a sufficient 
condition~\cite{Podewski}. It is based on the following generalization of 
the requirement that the complement to an arbitrary element 
must not be algebraic. 

\begin{definition}[Podewski \cite{Podewski}] A group $G$ is said to be 
\emph{ungebunden} if no system of fewer than $|G|$-many 
inequations in $G$ has 
precisely one solution. 
\end{definition}

The sufficient condition is that of being ungebunden.\footnote{Podewski 
stated in~\cite{Podewski} that any ungebunden algebraic system 
was topologizable. Palyutin, Seese, and Taimanov noticed that this 
is not true; they constructed a nontopologizable ungebunden 
ring~\cite{3}. However, 
Podewski's result is true for groups; moreover, a minor modification 
of the definition of the property of being ungebunden, which 
was suggested by Hesse and does not affect 
the case of groups, renders it true for 
arbitrary algebraic systems~\cite{Hesse}.} It is not necessary; Hesse proved 
the existence of a topologizable 
group of any uncountable cardinality $\lambda$ in which 
some countable system of inequations has precisely one solution~\cite{Hesse}. 
Hesse also suggested a generalization of the property of being ungebunden 
for arbitrary cardinals; namely, he called a group $\tau$-ungebunden if 
no system of $<\tau$ inequations in this group has precisely one solution. 
A group $G$ being ungebunden is equivalent to its being $|G|$-ungebunden, and 
the set $G\setminus\{1\}$ is additively algebraic in $G$ if and 
only if $G$ is $\aleph_0$-ungebunden. A natural generalization of this 
notion in terms of solutions to equations is as follows. 

\begin{definition}
Let $\tau$ be an arbitrary cardinal. 
We say that a subset $A$ of a group $G$ with identity element $1$  
is \emph{$\tau$-additively 
algebraic} if it can be represented as a union of 
fewer than $\tau$-many elementary algebraic sets. 
Arbitrary intersections of $\tau$-additively 
algebraic sets are called \emph{$\tau$-algebraic}. 
The $\aleph_0$-(additively) 
algebraic sets coincide with the (additively) algebraic sets. 
\end{definition} 

Recall that a topology  $\mathcal T$ is called a $P_\tau$-topology 
if any intersection of $<\tau$-many 
$\mathcal T$-open sets is $\mathcal T$-open. 

\begin{definition}
We say that $A$ is 
\emph{unconditionally $\tau$-closed} in $G$ if it is  closed in 
any Hausdorff group $P_\tau$-topology on $G$. 
Clearly, all unconditionally closed sets are unconditionally $\tau$-closed 
for any cardinal $\tau$, and unconditionally $\aleph_0$-closed sets 
coincide with unconditionally closed sets. 
By an 
\emph{$F_{<\tau}$-set} (in a given topology) 
we mean a union of $<\tau$-many closed sets. 
\end{definition}

It is easy to see that all $\tau$-algebraic sets are 
unconditionally $\tau$-closed; as mentioned above, 
the converse is not always true. In this paper, we prove, in 
particular, 
that all unconditionally $\tau$-closed sets are $\tau$-algebraic 
in any group whose center has index at most $\tau$.
 
\begin{definition}[Markov~\cite{Markov1946}] 
Let 
$m$ be a positive integer. By a \emph{multiplicative function} of 
$m$ arguments we mean an arbitrary element of the 
free product $\langle t_1\rangle * \dots * \langle t_m\rangle$, i.e., 
an irreducible 
word on the alphabet $\{t_1^{\pm1}, 
\dots, t_m^{\pm1}\}$ (the $t_i$ are treated as variables). 
The \emph{length} of the function is equal to 
the length of this word. Suppose that $G$ is a group and  
$g_1, \dots, g_m\in G$. The \emph{value} of a multiplicative function 
$\Phi=\prod_{i=1}^n t_{j_i}^{\varepsilon_i}$ of $m$ 
arguments (it is assumed that $j_i\le m$ for $i= 1, \dots, n$) 
at $g_1, \dots, g_m$ in $G$ is defined as 
$$
\Phi(g_1, \dots, g_m) =\prod_{i=1}^n g_{j_i}^{\varepsilon_i}.
$$
\end{definition}

Note that the set of all multiplicative functions 
is countable, and 
the set of all multiplicative functions of length 
$k$ with at most $n$ variables is finite for 
any positive integers $k$ and $n$. 

\begin{definition}
For a subset $A$ of a group $G$, 
$\mathop{\smash{\widetilde A}}\nolimits^{\,G}_\tau$ 
(or simply $\widetilde A$) denotes the \emph{$\tau$-algebraic closure} 
of $A$ in $G$, i.e., the intersection of all $\tau$-algebraic sets in 
$G$ containing $A$.  
\end{definition}

As 
usual, $[G]^{\le \tau}$ denotes the family of all subsets of $G$ of 
cardinality at most $\tau$. 
This is a partially ordered set under 
inclusion. If $\mathcal P$ is a property of groups such that any 
$X\in [G]^{\le \tau}$ is contained 
in a subgroup of cardinality $\le \tau$ with property $\mathcal P$, 
then set theorists would 
say that subgroups with property 
$\mathcal P$ form an \emph{unbounded set} in 
$[G]^{\le \tau}$. If, in addition,  
the union of any increasing chain of at most $\tau$ (precisely 
$\tau$) subgroups of cardinality $\le \tau$ with property $\mathcal 
P$ has property $\mathcal P$, then set theorists (we) say that 
subgroups with property $\mathcal P$ form a \emph{club} 
(\emph{$\tau$-club}) in $[G]^{\le 
\tau}$; \emph{club} is an abbreviation for 
\emph{closed unbounded}~(see, e.g., \cite[p.~916]{handbook}). Clearly, any 
\emph{club} is a $\tau$-\emph{club}. 

\begin{theorem}\label{thm1} 
Suppose 
that $\tau$ is a regular infinite cardinal, $G$ is a group, 
$\mathcal P$ is a property of groups such that subgroups with property 
$\mathcal P$ form a $\tau$-\emph{club} in 
$[G]^{\le \tau}$, and  
$A\subset G$ is not $\tau$-algebraic in $G$. 
Then the family of subgroups $G'\in [G]^{\le \tau}$ 
with property $\mathcal P$ 
in which $A\cap G'$ is not unconditionally $\tau$-closed is unbounded 
in $[G]^{\le \tau}$.
\end{theorem} 

\begin{proof}
For 
each $X\in [G]^{\le \tau}$, 
we fix  some subgroup $[[X]]$ of cardinality $\le \tau$ 
with property $\mathcal P$ containing $X$. 
For $x\in G$, we write 
$[[x]]$ instead of $[[\{x\}]]$. 

Let 
$X\in [G]^{\le \tau}$. 
Suppose that $A\subset G$, $1\in \widetilde A$, and 
$1\notin A$ (this implies that $A$ is not $\tau$-algebraic in $G$). 
Let us show that $A\cap G'$ is not unconditionally $\tau$-closed 
in some subgroup $G'$  
of $G$ with property $\mathcal P$ of cardinality $\le \tau$ which 
contains $X$. 
Take an arbitrary element $a_0\in A$. 
Let $\mathfrak M_1$ be the 
(finite) set of all finite sequences of the form 
$(\Phi, \underbrace {a_0,\dots, a_0}_{\text{$n$ times}})$, 
where $n\le 4$, 
$\Phi$ is a  multiplicative function of length $\le 6$ with  
$n+1$ arguments, and 
$$ 
\Phi(a_0,\dots, a_0, 1)\ne 1. 
$$ 
For any such sequence $(\Phi, a_0,\dots , a_0)$, 
we set 
$$
A_{\Phi,a_0, \dots, a_0}=
\{x\in G:\Phi(a_0,\dots, a_0, x)=1\}
$$
and 
$$
B_1= \bigcup_{(\Phi,a_0,\dots, a_0)
\in \mathfrak M_1} 
A_{\Phi,a_0,\dots, a_0}.
$$
The set $\mathfrak M_1$ is finite; hence  
$B_1$ is additively algebraic and, therefore, 
$\tau$-additively algebraic. 
We have $1\notin B_1$. On the other hand, $1\in \widetilde A$. Since 
$\widetilde A$ is the $\tau$-algebraic closure of $A$, we have 
$A\setminus B_1\ne \varnothing$. Take $x_1\in A\setminus B_1$. 
Note that $x_1\ne 1$ and $x_1\ne a_0$. Indeed, $x_1\ne 1$ 
because $x_1\in A\not\ni 1$, and $x_1\ne a_0$ because 
$\Phi=t_1t_2^{-1}$ 
is a multiplicative function (of length~2 with 
two arguments) for which $\Phi(a_0, 1) = a_0\ne 1$ 
(i.e., $(\Phi, a_0)\in \mathfrak M_1$)
but $\Phi(a_0, a_0) = a_0a_0^{-1} = 1$.

Let us represent the infinite cardinal $\tau$ (as usual, we assume that 
$\tau$ is the set of all ordinals of smaller cardinality) as the 
disjoint union  
of $\tau$-many of its subsets 
$$
\tau= \bigcup_{0<\alpha\in \tau}T_\alpha, 
$$
where each $T_\alpha$ is a set of ordinals smaller than $\tau$, 
$T_\alpha\cap T_\beta=\varnothing $ for $\alpha\ne \beta$, 
and $|T_\alpha|= \tau$ for each $\alpha>0$. 
The subgroup 
$G_1=[[\{a_0, x_1\}\cup X]]$ (it was defined at the beginning of the proof)  
has at most $\tau$ elements. We set $a_1=x_1$ and enumerate the 
elements of $G_1$ different from 1, $a_0$, and $a_1$ (if they exist) 
by ordinals from $T_1$ different from 0 and 1. 
We obtain $G_1=\{1,a_0, a_1\}\cup \{a_\beta: \beta\in T'_1\}$, 
where $T'_1$ is some (possibly, empty) subset of $T_1$ not containing 
0 and 1. 
The subgroup 
$G_1$ may be finite, but it necessarily 
contains $a_0$ and $a_1$. 

Suppose that $\alpha\in \tau$ and, for all $\gamma< \alpha$, 
we have defined increasing 
subgroups $G_\gamma\in [G]^{\le\tau}$ 
with property $\mathcal P$
whose elements are enumerated as $G_\gamma= \{1\}\cup 
\{a_\beta: \beta\le \gamma\}\cup\{a_\beta: \beta\in \bigcup_{0<\delta\le \gamma} 
T'_\delta\}$, where $T'_\delta\subset T_\delta\setminus \{\beta\le \gamma\}$, 
and elements $x_\gamma\in (A\cap G_\gamma)\setminus 
\bigcup_{\beta<\gamma} G_\beta$ such that, whatever a multiplicative 
function $\Phi$ of length $\le 3|\gamma|^2\cdot 2^{|\gamma|}$ 
with $\le 2|\gamma|^2\cdot 2^{|\gamma|}+1$ 
arguments 
(for infinite $\gamma$, any multiplicative function meets these requirements), 
the equality $\Phi(b_{\gamma_1}, \dots, b_{\gamma_n}, x_\gamma)= 1$ 
for some $\gamma_1, \dots, \gamma_n< \gamma$ and $b_\delta\in \{a_\delta, 
x_\delta\}$
implies  $\Phi(b_{\gamma_1}, \dots, b_{\gamma_n}, 1)= 1$. 
Let us define $G_\alpha$ and $x_\alpha$. 
Let $\mathfrak M_\alpha$ be the 
set of all finite sequences $(\Phi, b_{\gamma_1}, \dots, b_{\gamma_n})$, 
where $n$ is a positive integer not exceeding 
$2|\alpha|^2\cdot 2^{|\alpha|}$, 
$\Phi$ is a  multiplicative function of length 
$\le 3|\alpha|^2\cdot 2^{|\alpha|}$ with  
$n+1$ arguments, 
$\gamma_1, \dots, \gamma_n< \alpha$, $b_\delta\in \{a_\delta, x_\delta\}$ 
for $\delta=\gamma_1, \dots, \gamma_n$, and 
$$ 
\Phi(b_{\gamma_1},\dots b_{\gamma_n}, 1) \ne 1. 
$$ 
Note 
that the constraints on the length and number of arguments of $\Phi$ are 
only needed if $\alpha$ is finite; for infinite $\alpha$, they are  void, 
and  $\mathfrak M_\alpha$ is simply the 
set of all finite sequences $(\Phi, b_{\gamma_1}, \dots, b_{\gamma_n})$, 
where $n$ is a positive integer, 
$\Phi$ is an arbitrary  multiplicative function 
$n+1$ arguments, 
$\gamma_1, \dots, \gamma_n< \alpha$, $b_\delta\in \{a_\delta, x_\delta\}$ 
for $\delta=\gamma_1, \dots, \gamma_n$, and 
$\Phi(b_{\gamma_1},\dots b_{\gamma_n}, 1) \ne 1$. 

For any such sequence $(\Phi, b_{\gamma_1},\dots b_{\gamma_n})$, 
we set 
$$
A_{\Phi,b_{\gamma_1},\dots b_{\gamma_n}}=
\{x\in G:\Phi(b_{\gamma_1},\dots b_{\gamma_n}, x)=1\} 
$$
and 
$$
B_\alpha=\bigcup_{(\Phi,b_{\gamma_1},\dots b_{\gamma_n})
\in \mathfrak M_\alpha} 
A_{\Phi,b_{\gamma_1},\dots b_{\gamma_n}}.
$$
Each set $A_{\Phi,b_{\gamma_1},\dots b_{\gamma_n}}$ is elementary 
algebraic, and the cardinality of $\mathfrak M_\alpha$ is less than $\tau$ 
(it is finite if $\alpha$ is finite and equals $|\alpha|$ for 
infinite $\alpha$); 
therefore, $B_\alpha$ is $\tau$-additively algebraic. 
We have $1\notin B_\alpha$. On the other hand, $1\in \widetilde A$. Since 
$\widetilde A$ is the $\tau$-algebraic closure of $A$, we have 
$A\setminus B_\alpha\ne \varnothing$. Take $x_\alpha\in A\setminus B_\alpha$. 
Note that $x_\alpha\ne a_\beta$ for $\beta <\alpha$. 
Indeed, if $x_\alpha= a_\beta$ for some $\beta<\alpha$, 
then $\Phi=t_1t_2^{-1}$ 
is a multiplicative function (of length~2 with 
two arguments) for which $\Phi(a_\beta, 1) = a_\beta\ne 1$ 
(i.e., $(\Phi, a_\beta)\in \mathfrak M_\alpha$) 
but $\Phi(a_\beta, x_\alpha) = a_\beta x_\alpha^{-1} = 1$. 
If $a_\alpha$ has not yet been defined, we set $a_\alpha= x_\alpha$. 
We also set $G_\alpha = 
[[\bigcup_{\gamma<\alpha} G_\gamma\cup 
\{x_\alpha\}]]$ 
and enumerate the elements of $G_\alpha$ that 
have not been enumerated at the preceding steps by ordinals 
from $T_\alpha$, so that 
$G_\alpha= \{a_\beta: \beta\le \alpha\}\cup 
\{a_\beta: \beta\in T'_\alpha\}$, where 
$T'_\alpha\subset T_\alpha\setminus\{\beta\le \alpha\}$.

After $\tau$ steps, we obtain an increasing chain 
of 
subgroups $G_\alpha \in [G]^{\le \tau}$ with property $\mathcal 
P$ ($\alpha<\tau$) 
and an element $x_\alpha\in (A\cap G_\alpha)\setminus  
\{a_\beta:\beta<\alpha\}$ for each $\alpha$; moreover, 
$\bigcup_{\alpha\in \tau} G_\alpha=\{a_\beta: \beta\in \tau\}$. For 
$\beta\in \tau$, we set 
\begin{multline*}
H'_\beta=\{g^{-1}x^\varepsilon_\gamma g: 
\gamma\ge \beta,\ g=a_{\delta_1}\dots a_{\delta_n},\\ 
\text{where $n<\omega$, $\beta+n\le \gamma$, $\varepsilon=\pm1$, 
and $\delta_i<\gamma$}\}
\end{multline*}
and
\begin{multline*}
H_\beta=\bigcup 
\{H'_{\beta+\gamma_1}H'_{\beta+\gamma_2}\dots
H'_{\beta+\gamma_n}:\\ 
\text{$n\in \omega$ and $\gamma_1, \dots, \gamma_n$ are 
ordinals such that}\\ 
\text{$\beta+\gamma_i\in \tau$ and each $\gamma_i$ occurs
at most $2^{|\gamma_i|}$ 
times}\}.
\end{multline*}
Thus, 
$H_\beta$ is the union of products of  a certain form; namely, 
each product 
contains at most $2^n$ factors 
$H'_{\beta+n}$ with finite $n$ and arbitrarily many 
factors $H_{\beta+\gamma}$ with infinite $\gamma$ (provided that 
$\beta+\gamma<\tau$).

Note that $1\in H_\beta$ for each $\beta$ 
and $H_{\beta'}\subset H_{\beta''}$ if $\beta'<\beta''$. 

We 
claim that, if 
$\alpha<\beta\in \tau$, then 
$a_\alpha\notin H_\beta$. 
Indeed, suppose that, on the contrary, 
$a_\alpha \in H_\beta$ for $\beta> \alpha$. 
Then 
$a_\alpha = g_1^{-1}x_{\gamma_1}^{\varepsilon_1}
g_1\allowbreak
\dots \allowbreak 
g_n^{-1}x_{\gamma_n}^{\varepsilon_n}g_n$, where   
$\gamma_i \ge\beta$,  each factor 
$g_i^{-1}x_{\gamma_i}^{\varepsilon_i}g_i$ 
belongs to $H'_{\beta+\gamma'_i}$ 
(in particular, $\gamma_i\ge \beta+\gamma'_i$),  
every $\gamma'_i$ occurs at most $2^{|\gamma'_i|}$ 
times, 
$\varepsilon_i = \pm1$,   
and $g_i= a_{\delta_1(i)}\dots a_{\delta_{k_i}(i)}$, where 
$k_i$ is a positive integer such that $\beta+\gamma'_i+k_i\le \gamma_i$ 
and $\delta_j(i)< \gamma _i$ for $j\le k_i$. 
Suppose that $a_\alpha$ has no shorter (with smaller $n$)
representation in this form.  Let 
$\gamma^*= \max\{\gamma_1, \dots, \gamma_n\}$, 
and let $\{i_1, \dots, i_k\}= \{i\le n: \gamma_i= \gamma^*\}$. 
Note that $\beta+\gamma_i'\le \gamma^*$ for all $i\le n$; therefore, 
if $\gamma^*$ is finite, then 
 the assumption that 
every $\gamma'_i$ occurs at most $2^{|\gamma'_i|}$ times implies 
$n\le \gamma^*\cdot 2^{\gamma^*}$. Clearly, $k_i<\gamma^*$ for any $i$. 
Consider the obviously defined multiplicative function $\Phi$ 
of length $1+2(k_1+\dots+k_n)+n$ with $2 + k_1+\dots+k_n + n-k$ arguments 
for which 
\allowdisplaybreaks
\begin{align*}
&\Phi(a_\alpha,a_{\delta_1(1)}, \dots, a_{\delta_{k_1}(1)}, x_{\gamma_1},
a_{\delta_1(2)}, \dots, a_{\delta_{k_2}(2)}, x_{\gamma_2},\\ 
&\qquad\qquad
\dots, \\
&\qquad\qquad
a_{\delta_1(i_1-1)}, \dots, a_{\delta_{k_{i_1-1}}(i_1-1)}, x_{\gamma_{i_1-1}},
a_{\delta_1(i_1)}, \dots, a_{\delta_{k_{i_1}}(i_1)}, \\
&\qquad\qquad\qquad\qquad\qquad\qquad\qquad\qquad
a_{\delta_1(i_1+1)}, \dots, a_{\delta_{k_{i_1+1}}(i_1+1)}, x_{\gamma_{i_1+1}},\\
&\qquad\qquad
\dots, \\
&\qquad\qquad
a_{\delta_1(i_2-1)}, \dots, a_{\delta_{k_{i_2-1}}(i_2-1)}, x_{\gamma_{i_2-1}},
a_{\delta_1(i_2)}, \dots, a_{\delta_{k_{i_2}}(i_2)}, \\
&\qquad\qquad\qquad\qquad\qquad\qquad\qquad\qquad
a_{\delta_1(i_2+1)}, \dots, a_{\delta_{k_{i_2+1}}(i_2+1)}, x_{\gamma_{i_2+1}},\\
&\qquad\qquad
\dots, \\
&\qquad\qquad
a_{\delta_1(i_k-1)}, \dots, a_{\delta_{k_{i_k-1}}(i_k-1)}, x_{\gamma_{i_k-1}},
a_{\delta_1(i_k)}, \dots, a_{\delta_{k_{i_k}}(i_k)}, \\
&\qquad\qquad\qquad\qquad\qquad\qquad\qquad\qquad
a_{\delta_1(i_k+1)}, \dots, a_{\delta_{k_{i_k+1}}(i_k+1)}, x_{\gamma_{i_k+1}},\\
&\qquad\qquad
\dots, \\
&\qquad\qquad
a_{\delta_1(n)}, \dots, a_{\delta_{k_n}(n)}, x_{\gamma_{n}},
x_{\gamma^*})\\
&=a_\alpha^{-1}
a_{\delta_{k_1}(1)}^{-1} \dots a_{\delta_1(1)}^{-1} 
x_{\gamma_1}^{\varepsilon_1}
a_{\delta_1(1)} \dots a_{\delta_{k_1}(1)}\,
\dots \,
a_{\delta_{k_n}(n)}^{-1} \dots a_{\delta_1(n)}^{-1} 
x_{\gamma_n}^{\varepsilon_n}
a_{\delta_1(n)} \dots a_{\delta_{k_n}(n)}\\
&=a_\alpha^{-1} g_1^{-1}x_{\gamma_1}^{\varepsilon_1}
g_1\dots  
g_n^{-1}x_{\gamma_n}^{\varepsilon_n}g_n
=1.
\end{align*}
We have $\alpha<\beta\le \gamma^*$, 
$\delta_j(i)<\gamma^*$, 
$\gamma_i<\gamma^*$ for 
$i\notin\{i_1, \dots, i_k\}$, and 
$k_i< \gamma_i\le \gamma^*$ for all $i$. 
If $\gamma^*$ is finite, then, as mentioned 
above, $n\le \gamma^*\cdot 2^{\gamma^*}$ and $k_i<\gamma^*$; therefore, 
the length of $\Phi$ is no longer 
than $2(\gamma^*)^2\cdot 2^{\gamma^*} + \gamma^*\cdot 2^{\gamma^*} 
\le 3(\gamma^*)^2\cdot 2^{\gamma^*}$, 
and the number of its arguments is 
at most $(\gamma^*)^2\cdot 2^{\gamma^*}+\gamma^*\cdot 2^{\gamma^*}+1\le 
2(\gamma^*)^2\cdot 2^{\gamma^*}+1$.  
If $\gamma^*$ is infinite, then, whatever the length of 
$\Phi$, it is shorter than  
$3(|\gamma^*|)^2\cdot 2^{|\gamma^*|}$, and $\Phi$ has fewer than 
$2(|\gamma^*|)^2\cdot 2^{|\gamma^*|}+1$ arguments. 
It follows that 
\allowdisplaybreaks
\begin{align*}
&\Phi(a_\alpha,a_{\delta_1(1)}, \dots, a_{\delta_{k_1}(1)}, x_{\gamma_1},
a_{\delta_1(2)}, \dots, a_{\delta_{k_2}(2)}, x_{\gamma_2},\\ 
&\qquad\qquad
\dots, \\
&\qquad\qquad
a_{\delta_1(i_1-1)}, \dots, a_{\delta_{k_{i_1-1}}(i_1-1)}, x_{\gamma_{i_1-1}},
a_{\delta_1(i_1)}, \dots, a_{\delta_{k_{i_1}}(i_1)}, \\
&\qquad\qquad\qquad\qquad\qquad\qquad\qquad\qquad
a_{\delta_1(i_1+1)}, \dots, a_{\delta_{k_{i_1+1}}(i_1+1)}, x_{\gamma_{i_1+1}},\\
&\qquad\qquad
\dots, \\
&\qquad\qquad
a_{\delta_1(i_2-1)}, \dots, a_{\delta_{k_{i_2-1}}(i_2-1)}, x_{\gamma_{i_2-1}},
a_{\delta_1(i_2)}, \dots, a_{\delta_{k_{i_2}}(i_2)}, \\
&\qquad\qquad\qquad\qquad\qquad\qquad\qquad\qquad
a_{\delta_1(i_2+1)}, \dots, a_{\delta_{k_{i_2+1}}(i_2+1)}, x_{\gamma_{i_2+1}},\\
&\qquad\qquad
\dots, \\
&\qquad\qquad
a_{\delta_1(i_k-1)}, \dots, a_{\delta_{k_{i_k-1}}(i_k-1)}, x_{\gamma_{i_k-1}},
a_{\delta_1(i_k)}, \dots, a_{\delta_{k_{i_k}}(i_k)}, \\
&\qquad\qquad\qquad\qquad\qquad\qquad\qquad\qquad
a_{\delta_1(i_k+1)}, \dots, a_{\delta_{k_{i_k+1}}(i_k+1)}, x_{\gamma_{i_k+1}},\\
&\qquad\qquad
\dots, \\
&\qquad\qquad
a_{\delta_1(n)}, \dots, a_{\delta_{k_n}(n)}, x_{\gamma_{n}},
1)\\
&=a_\alpha^{-1}
a_{\delta_{k_1}(1)}^{-1} \dots a_{\delta_1(1)}^{-1} 
x_{\gamma_1}^{\varepsilon_1}
a_{\delta_1(1)} \dots a_{\delta_{k_1}(1)}\\
&\qquad\qquad\dots\\
&\qquad\qquad
a_{\delta_{k_{i_1-1}}(i_1-1)}^{-1}\dots a_{\delta_1(i_1-1)}^{-1} 
x_{\gamma_{i_1-1}}^{\varepsilon_{i_1-1}}
a_{\delta_1(i_1-1)} \dots a_{\delta_{k_{i_1-1}}(i_1-1)}\\
&\qquad\qquad
 a_{\delta_{k_{i_1+1}}(i_1+1)}^{-1} \dots a_{\delta_1(i_1+1)}^{-1}
x_{\gamma_{i_1+1}}^{\varepsilon_{i_1+1}}
a_{\delta_1(i_1+1)} \dots a_{\delta_{k_{i_1+1}}(i_1+1)}\\
&\qquad\qquad\dots\\
&\qquad\qquad
a_{\delta_{k_{i_k-1}}(i_k-1)}^{-1} \dots a_{\delta_1(i_k-1)}^{-1}
x_{\gamma_{i_k-1}}^{\varepsilon_{i_k-1}}
a_{\delta_1(i_k-1)} \dots a_{\delta_{k_{i_k-1}}(i_k-1)}\\
&\qquad\qquad
a_{\delta_{k_{i_k+1}}(i_k+1)}^{-1} \dots a_{\delta_1(i_k+1)}^{-1}
x_{\gamma_{i_k+1}}^{\varepsilon_{i_k+1}}
a_{\delta_1(i_k+1)} \dots a_{\delta_{k_{i_k+1}}(i_k+1)}\\
&\qquad\qquad\dots\\
&\qquad\qquad
a_{\delta_{k_n}(n)}^{-1} \dots a_{\delta_1(n)}^{-1}
x_{\gamma_n}^{\varepsilon_n}
a_{\delta_1(n)} \dots a_{\delta_{k_n}(n)}\\
&
=a_\alpha^{-1}
g_1^{-1}x_{\gamma_1}^{\varepsilon_1}g_1 
\dots  
g_{i_1-1}^{-1}x_{\gamma_{i_1-1}}^{\varepsilon_{i_1-1}}g_{i_1-1}
g_{i_1+1}^{-1}x_{\gamma_{i_1+1}}^{\varepsilon_{i_1+1}}g_{i_1+1}\\
&\qquad\qquad\qquad\qquad
\dots
g_{i_2-1}^{-1}x_{\gamma_{i_2-1}}^{\varepsilon_{i_2-1}}g_{i_2-1}
g_{i_2+1}^{-1}x_{\gamma_{i_2+1}}^{\varepsilon_{i_2+1}}g_{i_2+1}
\dots
g_{i_k-1}^{-1}x_{\gamma_{i_k-1}}^{\varepsilon_{i_k-1}}g_{i_k-1}\\
&\qquad\qquad\qquad\qquad\qquad\qquad
g_{i_k+1}^{-1}x_{\gamma_{i_k+1}}^{\varepsilon_{i_k+1}}g_{i_k+1}
\dots
g_n^{-1}x_{\gamma_n}^{\varepsilon_n}g_n =1, 
\end{align*}
i.e., 
\begin{multline*}
a_\alpha= 
g_1^{-1}x_{\gamma_1}g_1 
\dots 
g_{i_1-1}^{-1}x_{\gamma_{i_1-1}}g_{i_1-1}
g_{i_1+1}^{-1}x_{\gamma_{i_1+1}}g_{i_1+1}\\
\dots
g_{i_2-1}^{-1}x_{\gamma_{i_2-1}}g_{i_2-1}
g_{i_2+1}^{-1}x_{\gamma_{i_2+1}}g_{i_2+1}
\dots
g_{i_k-1}^{-1}x_{\gamma_{i_k-1}}g_{i_k-1}\\
g_{i_k+1}^{-1}x_{\gamma_{i_k+1}}g_{i_k+1}
\dots
g_n^{-1}x_{\gamma_n}g_n. 
\end{multline*}
We have obtained a shorter representation of $a_\alpha$ as 
an element of $H_\beta$. This contradiction proves the claim
that $a_\alpha\notin H_\beta$. 

Let $G'=\bigcup_{\alpha\in \tau} G_\alpha$. Then $G'$ is a 
subgroup 
of $G$ of cardinality $\le \tau$ with property $\mathcal P$ (because 
so are all $G_\alpha$ and the family of subgroups with property $\mathcal P$ 
forms a $\tau$-\emph{club} in $[G]^{\le\tau}$), 
$G'\supset X$, and the $H_\alpha$ with $\alpha\in \tau$ 
form a decreasing chain of subsets of $G'$ containing $1$. Moreover, 
$G'=\{1\}\cup \{a_\alpha: \alpha\in \tau\}$, 
and for 
each $a_\alpha$, there is a $\beta\in \tau$ such that 
$a_\alpha\notin H_\beta$. Therefore, 
$\bigcap_{\alpha\in \tau} H_\alpha=\{1\}$. Note that if 
$\alpha<\beta<\gamma\in \tau$, then $H_\alpha^{-1}= H_\alpha$, 
$H_\gamma^2\subset H_\beta$, and 
$a_\alpha^{-1}H_\gamma a_\alpha\subset H_\beta$. Finally, 
for any $a\in H_\alpha$, there exists a $\delta \in \tau$ for 
which $a H_\delta\subset H_\alpha$. Indeed, we have 
$a\in H'_{\alpha+\gamma_1}H'_{\alpha+\gamma_2}\dots
H'_{\alpha+\gamma_n}$ for some 
$n\in \omega$ and ordinals $\gamma_1, \dots, \gamma_n$ 
such that $\alpha+\gamma_i\in \tau$ and each $\gamma_i$ 
occurs at most $2^{|\gamma_i|}$ times. Clearly, for $\delta= \alpha+ 
\max\{\gamma_1, \dots, \gamma_n\}$, we have  
$a H_\delta\subset H_\alpha$. 

Thus, the sets $H_\beta$ with $\beta\in \tau$ form a base for some 
nondiscrete Hausdorff group topology $\mathcal T$ 
on $G'$. 
Since the intersection $H_\beta\cap A$ is nonempty 
(it contains $x_\beta$) for each 
$\beta$, it follows that 1 belongs to the closure of $A\cap G'$ 
in this topology, 
i.e, $A\cap G'$ 
is not unconditionally closed in $G'$.
Note that any intersection of 
less than $\tau$ sets $H_\alpha$ contains some $H_\beta$ (because these 
sets decrease and $\tau$ is regular), i.e., is open 
in $\mathcal T$. Hence any $F_{<\tau}$-set 
in 
$\mathcal T$ is closed in $\mathcal T$ (i.e., $\mathcal T$ is a 
$P_\tau$-topology), and $A\cap G'$ is not 
even unconditionally $\tau$-closed in $G'$. 

Now, let $B$ be an arbitrary non-$\tau$-algebraic set in $G$, 
i.e., a set such that 
$\widetilde B\ne B$. Take $b\in \widetilde B\setminus 
B$. We have $1\notin b^{-1} B$. On the other hand, $1\in 
b^{-1}\widetilde B$. It is easy to see that  
$b^{-1} {\widetilde B}= \widetilde{b^{-1}B}$ 
(the proof is similar to that of 
Lemma~12 from \cite{Markov1946}); hence $1\in \widetilde{b^{-1}B}\setminus 
b^{-1}B$. As was shown above, this implies 
the existence of 
a subgroup $G'\in [G]^{\le \tau}$ with 
property $\mathcal P$ such that it contains 
$b$ (the one-point set $\{b\}$ 
plays the role of $X$) and $1$ belongs to the closure of 
$b^{-1}B\cap G'$ in 
some group $P_\tau$-topology $\mathcal T$ on $G'$. Since $b\in G'$, 
it follows that $b(b^{-1})B\cap G'=B\cap G'$, and the point $b$ 
belongs to the 
closure of $B\cap G'$ in $\mathcal T$, i.e., $B\cap G'$ 
is not closed in $G'$ with the topology $\mathcal T$.
\end{proof}

\begin{remark}\label{rem1}
Any 
algebraic 
set in the group $G$ is unconditionally closed in 
$G$~\cite[Theorem~1]{Markov1946}. Therefore, any $\tau$-algebraic 
set in $G$ is an intersection of $F_{<\tau}$-sets 
in any Hausdorff group topology on $G$ and hence 
unconditionally $\tau$-closed.
\end{remark} 

\begin{remark}\label{rem2}
The regularity of $\tau$ is needed only to show that $A$ is not 
unconditionally $\tau$-closed. Without this assumption, we obtain the 
following statement: \emph{Suppose 
that $\tau$ is an infinite cardinal, $G$ is a group, 
$\mathcal P$ is a property of groups such that subgroups with property 
$\mathcal P$ form a $\tau$-\emph{club} in 
$[G]^{\le \tau}$, and  
$A\subset G$ is not $\tau$-algebraic in $G$. 
Then the family of subgroups $G'\in [G]^{\le \tau}$ with property $\mathcal P$ 
in which $A\cap G'$ is not unconditionally closed is unbounded 
in $[G]^{\le \tau}$.}
\end{remark}

\begin{definition}[D. Dikranjan and D. Shakhmatov]
A 
subgroup $H$ of a group $G$ is said to be \emph{supernormal} 
in $G$ if, for any $g\in G$, there exists an $h\in H$  such that 
$g^{-1}xg=h^{-1}xh$ for all $x\in H$. 
\end{definition}

\begin{corollary}\label{cor1}
Suppose 
that $\tau$ is a regular infinite cardinal and $G$ is a group 
such that $[G]^{\le \tau}$ contains a 
$\tau$-\emph{club} consisting of supernormal subgroups. 
Then any set 
$A\subset G$ is $\tau$-algebraic in $G$ if and only if 
$A$ is unconditionally $\tau$-closed in $G$.
\end{corollary}

\begin{proof}
Suppose that $A\subset G$ and $A$ is not $\tau$-algebraic in $G$. 
Let $\mathcal P$ be the property of belonging to the given 
$\tau$-\emph{club} of supernormal subgroups. 
By Theorem~1, there exists a supernormal subgroup 
$H$ in $G$ in which 
$A\cap H$ is not unconditionally $\tau$-closed, 
i.e., there exists  a Hausdorff 
group $P_\tau$-topology $\mathcal T$ on $H$ in which $A\cap H$ is not closed. 
Let $\mathcal B$ be any  neighborhood base at the identity element 
of $\mathcal T$. Since $H$ is supernormal, it follows that $\mathcal B$ 
is also a base for some group topology on $G$; clearly, this topology is 
Hausdorff and $P_\tau$, and $A$ is not closed in it. 
Thus, any unconditionally $\tau$-closed subset of $G$ is $\tau$-algebraic 
in $G$. The converse implication holds by Remark~\ref{rem1}.
\end{proof}

\begin{corollary}[see also \cite{sipa}]\label{cor2}
A 
subset of a direct product $G$ of an Abelian group $A$ and a direct product 
$B=\prod_{\alpha \in I}B_\alpha$ of countable groups  
is unconditionally closed in $G$ if and only if 
it is algebraic in $G$. 
\end{corollary}

\begin{proof}
We say that $B'$  is a countable subproduct in $B$ if 
$B=\prod_{\alpha \in I}B'_\alpha$, where $B'_\alpha=B_\alpha$ for countably 
many indices $\alpha$ and $B'_\alpha=\{1_\alpha\}$ (here $1_\alpha$ is the 
identity element of $B_\alpha$) for the remaining $\alpha$. Clearly, all 
(countable) subgroups $A'\times B'$ of $G$, where $A'$ is a countable subgroup 
of $A$ and $B'$ is a countable subproduct in $B$, are supernormal in $G$, and 
they form a \emph{club} in $[G]^{\le \omega}$. It remains to apply 
Corollary~\ref{cor1}. 
\end{proof}

\begin{corollary}\label{cor3}
If 
$\tau$ is an infinite regular cardinal and 
$G$ is a group with center of index $\le \tau$, then 
a set 
$A\subset G$ is $\tau$-algebraic in $G$ if and only if 
$A$ is unconditionally $\tau$-closed in $G$.
\end{corollary}

\begin{proof}
Let $Z$ be the center of $G$. The condition $|G:Z|\le \tau$ means that 
there exists a set $X\subset G$ of cardinality $\le\tau$ such that, 
for any $g\in G$, there exists an $x\in X$ and a $z\in Z$ for which $g=xz$. 
Any subgroup $H\subset G$ containing $X$ is supernormal in $G$. Indeed, 
if $g=xz$, where  $x\in X(\subset H)$ and $z\in Z$, 
then, for any $h\in H$, we have 
$g^{-1}hg=z^{-1}x^{-1}hxz=x^{-1}hx$. It remains to note that 
the family of all subgroups $H\in [G]^{\le \tau}$ containing $X$ forms a 
\emph{club} in $[G]^{\le \tau}$ and apply Corollary~\ref{cor1}.
\end{proof}

Since the center of any group $G$ has index $\le |G|$, we obtain 
the following corollary (see also Remark~\ref{rem2}).

\begin{corollary}[Podewski \cite{Podewski}]
If the set $G\setminus \{1\}$ 
is not $|G|$-algebraic in $G$ 
\textup(i.e., $G$ is ungebunden\textup), 
then $G\setminus \{1\}$ is not unconditionally closed 
\textup(i.e., $G$ admits a nondiscrete Hausdorff group topology\textup). 
\end{corollary}

\begin{corollary}\label{cor5}
An infinite  
group $G$ of regular cardinality 
is ungebunden if and only if it 
admits a nondiscrete Hausdorff group $P_{|G|}$-topology. 
\end{corollary}

\subsection*{Acknowledgments}
The 
author is very grateful to Dikran Dikranjan and Dmitri Shakhmatov, who noticed 
a mistake in the original version of this paper, and to the 
referee for valuable suggestions.


\begin{thebibliography}{10} 

\bibitem{handbook}
J. E. Baumgartner, 
``Applications of the Proper Forcing Axiom,''
in \textit{Handbook of Set-Theoretic Topology}, 
Ed. by K.~Kunen and J.~E.~Vaughan (North-Holland, Amsterdam, 1984), 
pp.~913--959. 

\bibitem{Hesse}
G. Hesse,  \textit{Zur Topologisierbarkeit von Gruppen}, 
Dissertation (Univ. Hannover, Hannover, 1979).

\bibitem{Markov1944} 
A. A. Markov. 
``On unconditionally closed sets,'' 
Dokl. Akad. Nauk SSSR 
\textbf{44} (5), 196--197 (1944); English translation: 
A. Markoff, 
``On unconditionally closed sets,'' 
C. R. (Doklady) Acad. Sci. URSS (N.S.) 
\textbf{44}, 180--181  (1944).
 
\bibitem{Markov1945} 
A. A. Markov, 
``On free topological groups,''
Izv. Akad. Nauk SSSR, Ser. Mat. \textbf{9} (1), 3--64 (1945); 
 English translation: 
``Three papers on topological groups:  
I. On the existence of periodic  connected topological groups.  
II. On free topological groups.  
III. On  unconditionally closed sets,''
Amer. Math. Soc. Transl. \textbf{30} (1950).
 
\bibitem{Markov1946} 
A. A. Markov, 
``On unconditionally closed sets,'' 
Mat. Sb. \textbf{18} (1), 3--26 (1946); 
 English translation: 
``Three papers on topological groups:  
I. On the existence of periodic  connected topological groups.  
II. On free topological groups.  
III. On  unconditionally closed sets,''
Amer. Math. Soc. Transl. \textbf{30} (1950).

\bibitem{3}
E. A. Paljutin, D. G. Seese, and A. D. Taimanov,
``A remark on the topologization of algebraic structures,''
Rev. Roumaine Math. Pures Appl. \textbf{26} (4), 617--618  (1981).

\bibitem{Podewski} 
K.-P. Podewski, 
``Topologisierung algebraischer Strukturen,''
Rev. Roumaine Math. Pures Appl. \textbf{22} (9), 1283--1290 (1977).


\bibitem{Sipa} 
O. V. Sipacheva, 
``Consistent solution of Markov's problem about algebraic sets,'' 
ArXiv:math.GR/0605558. 

\bibitem{sipa}
O. V. Sipacheva, ``A class of groups in which all 
unconditionally closed sets are algebraic,'' 
ArXiv:math.GR/0610430.


\end{thebibliography}
\end{document}